\documentclass[11pt]{article}
\usepackage{amsfonts}
\usepackage{latexsym}
\usepackage{amsmath}
\usepackage{amssymb}
\usepackage{amsthm}
\usepackage[affil-it]{authblk}
\usepackage[english]{babel}
\usepackage{}
\makeatletter \oddsidemargin -.1in \evensidemargin -.1in
\newcommand{\singlespacing}{\let\CS=\@currsize\renewcommand{\baselinestretch}{1}\tiny \CS}
\newcommand{\doublespacing}{\let\CS=\@currsize\renewcommand{\baselinestretch}{1.35}\tiny \CS}
\date{}
\newcommand{\be}{\begin{equation}}
\newcommand{\ee}{\end{equation}}
\newcommand{\bea}{\begin{eqnarray}}
\newcommand{\eea}{\end{eqnarray}}
\newcommand{\nn}{\nonumber}
\newcommand{\bee}{\begin{eqnarray*}}
\newcommand{\eee}{\end{eqnarray*}}

\newtheorem{Theorem}{Theorem}

\newtheorem{t1}{Theorem}[section]
\newtheorem{d1}{Definition}[section]

\newtheorem{l1}{Lemma}[section]
\newtheorem{p1}{Proposition}[section]

\begin{document}
\title{Controlled $K$-frames in Hilbert $C^*$-modules}
\author{Ekta Rajput{\dag}~, N. K.
Sahu\thanks{corresponding author}
\\
\small \it{{\dag}{*}Dhirubhai Ambani Institute of Information and
Communication Technology, Gandhinagar, India}}

\date{}
\maketitle \setlength{\baselineskip}{18pt} \noindent 
\begin{abstract}
Controlled frames have been the subject of interest because of its
ability to improve the numerical efficiency of iterative
algorithms for inverting the frame operator. In this paper, we
introduce the notion of controlled $K$-frame in Hilbert
$C^{*}$-modules. We establish the equivalent condition for
controlled $K$-frame. We investigate some operator theoretic
characterizations of controlled $K$-frames and controlled Bessel
sequences. Moreover we establish the relationship between the
$K$-frames and controlled $K$-frames. We also investigate the
invariance of a $C$-controlled $K$-frame under a suitable map $T$.
At the end we prove a perturbation result for controlled
$K$-frame.
\end{abstract}



\section{Introduction}
Frames a more flexible substitutes of bases in Hilbert spaces were
first proposed by Duffin and Schaeffer \cite{Duffin} in 1952 while
studying nonharmonic Fourier series. Daubechies, Grossmann and
Meyer \cite{Daubechies} reintroduced and developed the theory of
frames in 1986. Due to their rich structure the subject drew the
attention of many mathematician,
 physicists and engineers because of its applicability in signal processing \cite{Ferreira}, image processing \cite{Cand}, coding and communications \cite{Strohmer},
 sampling \cite{Eldar, Eldar Y.}, numerical analysis, filter theory \cite{olcskei}. Now a days it is used in compressive sensing, data analysis and other areas.
 In general frames can be viewed as a redundant representation of  basis. Due to its redundancy it becomes more applicable not only in theoretical point of view but
 also in various kinds of applications. \\
Hilbert $C^*$-modules are generalizations of Hilbert spaces by
allowing the inner product to take values in a $C^*$-algebra
rather than in the field of real or complex numbers. They were
introduced and investigated initially by Kaplansky
\cite{Kaplansky}. Frank and Larson \cite{Frank} defined the
concept of standard frames in finitely or countably generated
Hilbert $C^*$ -modules over a unital $C^*$ -algebra. For more
details of frames in Hilbert $C^{*}$-modules one may refer to
Doctoral Dissertation \cite{Jing}, Han et al. \cite{Han1} and Han
et al. \cite{Han2}. In 2012, L. Gavruta \cite{Gvruta} introduced
the notion of $K$-frames in Hilbert space to study the atomic
systems with respect to a bounded linear operator $K$. Controlled
frames in Hilbert spaces have been introduced by P. Balazs
\cite{Balazs} to improve the numerical efficiency of iterative
algorithms for inverting the frame operator. Rahimi \cite{Rahimi}
defined the concept of controlled $K$-frames in Hilbert spaces and
showed that controlled $K$-frames are equivalent to $K$-frames due
to which the controlled operator C can be used as preconditions in
applications. In \cite{Najati A.}, Najati et al. introduced the
concepts of atomic system for operators and $K$-frames in Hilbert
$C^*$-modules. Controlled frames in Hilbert $C^*$-modules were
introduced by Rashidi and Rahimi \cite{Rashidi}, and the authors
showed that they share many useful properties with their
corresponding notions in a Hilbert space. Motivated by the above
literature, we introduce the notion of a controlled $K$-frame in
Hilbert $C^*$-modules.
\section{Preliminaries}
In this section we give some basic definitions related to Hilbert
$C^{*}$-modules, frames, $K$-frames, Controlled frames in Hilbert
$C^*$-modules. Hilbert $C^*$-modules are generalization of Hilbert
spaces by allowing the inner product to take values in
$C^*$-algebra rather than $\mathbb{R}$ or $\mathbb{C}$.
\begin{d1}\rm
Let $\mathcal{A}$ be a $C^*$-algebra. An \textit{inner product $\mathcal{A}$-module} is a complex vector space $\mathcal{H}$ such that \\
(i) $\mathcal{H}$ is a right $\mathcal{A}$-module i.e there is a bilinear map
\begin{center}
    $\mathcal{H} \times \mathcal{A} \to \mathcal{A} \colon (x,a) \to x \cdot a$
\end{center}
satisfying $(x \cdot a) \cdot b = x \cdot (ab)$ and $(\lambda x) \cdot a = x \cdot (\lambda a)$, and $x \cdot 1 = x$ where $\mathcal{A}$ has a unit 1. \\
(ii) There is a map
$\mathcal{H} \times \mathcal{H} \to \mathcal{A} \colon (x,y) \to \langle x,y \rangle$ satisfying\\
1. $\langle x,x \rangle \geq 0$ \\
2. $\langle x,y \rangle^* = \langle y,x \rangle$\\
3. $\langle ax, y \rangle = a \langle x, y \rangle$\\
4. $\langle x+y, z \rangle = \langle x, z \rangle + \langle y,z \rangle$\\
5. $ \langle x,x \rangle = 0$ if and only if $x = 0$ (for every
$x, y, z \in \mathcal{H}$, $a \in \mathcal{A}$).
\end{d1}
\begin{d1}\rm  A \textit{Hilbert $C^*$-module} over $\mathcal{A}~$ is an inner product $\mathcal{A}$-module with the property that $(\mathcal{H},\| \cdot \|_\mathcal{H})$ is
complete with respect to the norm $\|x\| = \| \langle x, x \rangle
\|_{\mathcal{A}}^{\frac{1}{2}}$, where $\|.\|_{\mathcal{A}}$
denotes the norm on $\mathcal{A}$.
\end{d1}
Let $\mathcal{A}$ be a $C^*$-algebra and consider
\begin{center}
    $l^2 \mathcal{(A)} = \{\{a_j\} \subseteq \mathcal{A} : \sum_{j\in \mathbb{J}}a_{j}a_{j}^{*} \}.$
\end{center}
where the sum converges in norm in $\mathcal{A}$. It is easy to
see that $l^2 \mathcal{(A)}$ is a Hilbert $C^*$-module with
pointwise operations and the inner product defined as \bee \langle
\{a_j\}, \{b_j\} \rangle = \sum_{j\in \mathbb{J}} a_{j}b_{j}^{*},
~~~ \{a_j\}, \{b_j\} \in l^2 \mathcal{(A)} \eee and \bee
\|\{a_j\}\| = \sqrt{\|\sum_{j\in \mathbb{J}}  a_j a_{j}^{*}\|}.
\eee

\begin{d1}\rm (\cite{Jing})
Let $\mathcal{A}$ be a unital $C^*$-algebra and ${j \in
\mathbb{J}}$ be a finite or countable index set. A sequence
$\{\psi_j\}_{j \in \mathbb{J}}$ of elements in a Hilbert
$\mathcal{A}$ -module $\mathcal{H}$ is said to be a frame if there
exist two constants $C, D > 0$ such that \bea
    C \langle f,f \rangle \leq \sum_{j \in \mathbb{J}} \langle f,\psi_j \rangle \langle \psi_j,f \rangle \leq D\langle f,f \rangle,
    \forall f \in \mathcal{H}.
\eea
\end{d1}
The frame $\{\psi_j\}_{j \in \mathbb{J}}$ is said to be a tight
frame if $C = D$, and is said to be Parseval or a normalized tight
frame if $C = D = 1$.\\
Suppose that $\{\psi_j\}_{j \in \mathbb{J}}$ is a frame of a
finitely or countably generated Hilbert $C^{*}$-module
$\mathcal{H}$ over a unital $C^*$-algebra $\mathcal{A}$. The
operator $T \colon \mathcal{H} \to l^2\mathcal{(A)}$ defined by
\begin{center}
    $Tf=\{\langle f,\psi_j \rangle \}_{j\in \mathbb{J}}$
\end{center}
is called the \textit{analysis operator}.\\
The adjoint operator $T^* \colon l^2 \mathcal{(A)} \to \mathcal{H}$ is given by
\bee
    T^*\{c_j\}_{j \in \mathbb{J}}=\sum_{j \in \mathbb{J}}{c_j}{\psi_j}
\eee
$T^*$ is called \textit{pre-frame operator or the synthesis operator}.\\
By composing T and $T^*$, we obtain the \textit{frame operator} $S
\colon \mathcal{H}  \to \mathcal{H} $ \bea\label{eqn2}
    Sf=T^*Tf=\sum_{j \in \mathbb{J}}\langle f,\psi_j\rangle {\psi_j}.
\eea
\begin{d1}\rm \cite{Najati A.}
  A sequence $\{\psi_j\}_{j \in \mathbb{J}}$ of elements in a Hilbert  $\mathcal{A}$-module $\mathcal{H}$ is said to be a $K$-frame $(K \in L(\mathcal{H}))$ if there exist
  constants $C, D > 0$ such that
\bea
    C\langle K^*f, K^*f\rangle \leq \sum_{j \in \mathbb{J}} \langle f,\psi_j \rangle \langle \psi_j,f \rangle \leq D \langle f,f \rangle, \forall f \in \mathcal{H}.
\eea
\end{d1}
\begin{d1}\rm \cite{Rashidi}
Let $\mathcal{H} $ be a Hilbert $C^*$-module and $C \in
GL\mathcal{(H)}$. A frame controlled by the operator $C$ or $C$
-controlled frame in Hilbert $C^*$-module $\mathcal{H} $ is a
family of vectors $\{\psi_j\}_{j \in \mathbb{J}}$, such that there
exist two constants $A, B > 0$ satisfying \bee
    A \langle f, f \rangle \leq \sum_{j \in \mathbb{J}} \langle f,\psi_j \rangle  \langle C\psi_j,f \rangle \leq B \langle f,f \rangle, \forall f \in \mathcal{H}.
\eee Likewise, $\{\psi_j\}_{j \in \mathbb{J}}$ is called a $C$-controlled Bessel sequence with bound $B$, if there exists $B> 0
$ such that \bee
    \sum_{j \in \mathbb{J}} \langle f,\psi_j \rangle  \langle C\psi_j,f \rangle \leq B \langle f,f \rangle, \forall f \in
    \mathcal{H},
\eee
where the sum in the above inequalities converges in norm.\\
If $A = B$, we call $\{\psi_j\}_{j \in \mathbb{J}}$ as $C$-controlled tight frame, and if $A = B = 1$ it is called a
$C$-controlled Parseval frame.
\end{d1}

\section{Controlled operator frames }
For the rest of the paper we assume that $\mathcal{H}$ is a
Hilbert $C^{*}$-module over unital $C^{*}$-algebra $\mathcal{A}$
with $\mathcal{A}$-valued inner product $\langle .,.\rangle$ and
norm $\|.\|$. $L(\mathcal{H})$ denotes the set of all adjointable
operators on Hilbert $C^{*}$-module $\mathcal{H}$, and
$GL^{+}(\mathcal{H})$ indicates the set of all bounded linear
positive invertible operators on $\mathcal{H}$ with bounded
inverse.
 We define below the controlled operator frame or
$C$-controlled $K$-frame on a Hilbert $C^{*}$-module
$\mathcal{H}$.
\begin{d1}\rm Let $\mathcal{H}$ be a Hilbert $\mathcal{A}$-module over a unital $C^*$-algebra, $C \in GL^{+}(\mathcal{H})$ and $K \in L(\mathcal{H})$.
A sequence $\{\psi_j\}_{j \in \mathbb{J}}$ in $\mathcal{H}$ is
said to be a $C$-controlled $K$-frame if there exist two constants
$0<A\leq B< \infty$ such that \bea \label{e1.1}
    A \langle C^\frac{1}{2}K^*f, C^\frac{1}{2}K^*f \rangle \leq \sum_{j \in \mathbb{J}} \langle f,\psi_j \rangle  \langle C\psi_j,f \rangle \leq B \langle f,f \rangle, \forall f
    \in \mathcal{H}.
\eea
\end{d1}
\noindent If $C = I$, the $C$-controlled $K$-frame $\{\psi_j\}_{j
\in \mathbb{J}}$ is simply $K$-frame in $\mathcal{H}$ which was
discussed in \cite{Najati A.}. The sequence $\{\psi_j\}_{j \in
\mathbb{J}}$ is called a $C$-controlled Bessel sequence with bound
$B$, if there exists $B>0$ such that \bea \label{e1.2}
    \sum_{j \in \mathbb{J}} \langle f,\psi_j \rangle  \langle C\psi_j,f \rangle \leq B \langle f,f \rangle, \forall f \in
    \mathcal{H},
\eea
where the sum in the above inequalities converges in norm.\\
If $A = B$, we call this $C$-controlled $K$-frame a tight $C$-controlled $K$-frame, and if $A = B = 1$ it is called a Parseval $C$-controlled $K$-frame.\\
Let $\{\psi_j\}_{j \in \mathbb{J}}$ be a $C$-controlled Bessel sequence for Hilbert module $\mathcal{H}$ over $\mathcal{A}$.\\
The operator $T \colon {H} \to l^2(\mathcal{A})$ defined by
\bea \label{e1.3}
    Tf=\{\langle f,\psi_j \rangle\}_{j \in \mathbb{J}},    f \in \mathcal{H}
\eea is called the \textit{analysis operator}. The adjoint
operator $T^* \colon l^2(\mathcal{A}) \to \mathcal{H} $  given by
\bea \label{e1.4}
    T^*(\{c_j\})_{j \in \mathbb{J}}=\sum_{j \in \mathbb{J}}{c_j}C{\psi_j}
\eea is called \textit{pre-frame operator or the synthesis
operator}. By composing $T$ and $T^*$, we obtain the
\textit{C-controlled frame operator} $S_C \colon \mathcal{H} \to
\mathcal{H} $ as \bea \label{e1.5}
    S_Cf=T^*Tf=\sum_{j \in \mathbb{J}} \langle f,\psi_j \rangle
    C{\psi_j}.
\eea

We quote the following results from the literature that will be
used in our work.

\begin{l1} \label{l3.1} \rm \cite{Arambaic}  Let $\mathcal{A}$ be a $C^*$-algebra. Let $U$ and $V$ be two Hilbert $\mathcal{A}$-modules and $T\in End_{\mathcal{A}}^{*}(U,V)$.
Then the following statements are equivalent:
\begin{enumerate}
    \item $T$ is surjective.
    \item $T^*$ is bounded below with respect to norm i.e there exists $m > 0$ such that $\|T^{*}f\| \geq m \|f\|$ for all $f \in U$.
    \item  $T^*$ is bounded below with respect to inner product i.e there exists $m > 0$ such that $\langle T^{*}f,T^{*}f \rangle \geq m \langle f,f \rangle$ for all $f \in U$.
\end{enumerate}
\end{l1}
\begin{l1} \label{l3.2} \rm \cite{Paschke} Let $U$ and $V$ be Hilbert $\mathcal{A}$-modules over a $C^*$-algebra $\mathcal{A}$ and let $T:U \to V$ be a linear map.
Then the following conditions are equivalent:
\begin{enumerate}
    \item The operator $T$ is bounded and $\mathcal {A}$-linear.
    \item  There exists $k \geq 0$ such that $\langle Tx,Tx \rangle \leq k \langle x,x\rangle$ holds for all $x \in U$.

\end{enumerate}
\end{l1}
\begin{t1} \label{t3.1}\rm \cite{Fang} Let $E$, $F$ and $G$ be Hilbert $\mathcal{A}$-modules over a $C^*$-algebra $\mathcal{A}$. Let $T \in L(E,F)$ and $T^{'} \in L(G,F)$
 with $\overline{R(T^*)}$ be orthogonally complemented. Then the following statements are equivalent:
    \begin{enumerate}
        \item $T^{'}T^{'*} \leq \lambda TT^*$ for some $\lambda > 0$;
        \item There exists $\mu > 0$ such that $\|T^{'*}z\| \leq \mu \|T^*z\|$ for all $z \in F;$
        \item There exists $D \in L(G,E)$ such that $T^{'}=TD$, that is the equation $TX=T^{'}$ has a solution;
        \item $R(T^{'}) \subseteq R(T)$.
    \end{enumerate}
\end{t1}
For the rest of the paper we indicate that $S_{C}$ stands for the
controlled frame operator as we have defined in (\ref{e1.5}), and
$S$ stands for the classical frame operator in Hilbert
$C^{*}$-module $\mathcal{H}$ as defined in (\ref{eqn2}).
\begin{l1} \label{l3.5} \rm Let $C \in GL^{+}(\mathcal{H})$, $KC = CK$ and $R(C^{\frac{1}{2}}) \subseteq R(K^*C^{\frac{1}{2}})$ with $\overline{R((C^{\frac{1}{2}})^*) }$ is
 orthogonally complemented. Then $\| C^{\frac{1}{2}}f\|^2 \leq  \lambda^{'} \|K^*C^{\frac{1}{2}}f\|^2$ for some $\lambda^{'} > 0$.
\end{l1}
\begin{proof} Suppose $R(C^{\frac{1}{2}}) \subseteq R(K^*C^{\frac{1}{2}})$ with $\overline{R((C^{\frac{1}{2}})^*) }$ orthogonally complemented.
Then by using Theorem \ref{t3.1}, there exist some $\lambda^{'} >
0$ such that
    \bee
    (C^{\frac{1}{2}}){(C^{\frac{1}{2}})^*} \leq \lambda^{'} (K^*C^{\frac{1}{2}})
    (K^*C^{\frac{1}{2}})^*.
    \eee
   This implies that $\langle (C^{\frac{1}{2}}){(C^{\frac{1}{2}})^*}f,f \rangle \leq \lambda^{'}\langle (K^*C^{\frac{1}{2}}) (K^*C^{\frac{1}{2}})^*f,f\rangle$.\\
    Now by taking norm on both sides, we get
    \bee
    \| C^{\frac{1}{2}}f\|^2 \leq  \lambda^{'}
    \|K^*C^{\frac{1}{2}}f\|^2.
    \eee
\end{proof}
In the following theorem, we establish an equivalence condition
for $C$-controlled $K$-frame in a Hilbert $C^{*}$-module
$\mathcal{H}$.
\begin{t1} \label{t3.2} \rm Let $\mathcal{H}$ be a finitely or countably generated Hilbert $\mathcal{A}$ -module over a unital $C^*$-algebra
$\mathcal{A}$, $\{\psi_j\}_{j \in \mathbb{J}} \subset \mathcal{H}$
be a sequence, $C \in GL^{+}(\mathcal{H})$, $K \in
L(\mathcal{H})$, $KC=CK$ and $R(C^{\frac{1}{2}}) \subseteq
R(K^*C^{\frac{1}{2}})$ with $\overline{R((C^{\frac{1}{2}})^*) }$
be orthogonally complemented. Then $\{\psi_j\}_{j \in \mathbb{J}}$
is a $C$-controlled $K$-frame in Hilbert $C^*$-module if and only
if there exist constants $0< A\leq B< \infty$  such that
\bea\label{eqn1}
    A\| C^{\frac{1}{2}}K^{*}f\|^2  \leq \|\sum_{j \in \mathbb{J}} \langle f,\psi_j \rangle  \langle C\psi_j,f \rangle \| \leq B \| f \|^2,~ \forall f \in \mathcal{H}.
\eea
\end{t1}
\begin{proof}
($\implies$)  Obvious.\\
Now we assume that there exist constants $0<A,B< \infty $ such
that \bee
    A\| C^{\frac{1}{2}}K^{*}f\|^2  \leq \|\sum_{j \in \mathbb{J}} \langle f,\psi_j \rangle  \langle C\psi_j,f \rangle \| \leq B \| f
    \|^2,~~\forall f\in \mathcal{H}.
\eee We prove that $\{\psi_j\}_{j \in \mathbb{J}}$ is a
$C$-controlled $K$-frame for  Hilbert $C^*$-module $\mathcal{H}$.
As $S$ and $C$ are both positive operator, they are self adjoint.
Thus we have \bea \label{e1.6} A\| C^{\frac{1}{2}}K^{*}f\|^2
&\leq& \|\sum_{j \in \mathbb{J}} \langle f,\psi_j \rangle  \langle
C\psi_j,f \rangle \|\nn\\
&=& \|\langle S_{C}f, f\rangle\| =\|\langle
CSf,f\rangle\|=\|\big\langle (CS)^{\frac{1}{2}}f,
(CS)^{\frac{1}{2}}f \big\rangle\|, {\rm~as~} S_{C}=CS\nn\\
&=& \|(CS)^{\frac{1}{2}}f\|^{2}.\eea Since $R(C^{\frac{1}{2}})
\subseteq R(K^*C^{\frac{1}{2}})$ with
$\overline{R((C^{\frac{1}{2}})^*) }$ is orthogonally complemented,
then using Lemma \ref{l3.5},
 there exist some $\lambda^{'} > 0$ such that
\bee \| C^{\frac{1}{2}}f\|^2 \leq  \lambda^{'}
\|K^*C^{\frac{1}{2}}f\|^2. \eee Multiplying both side by $A$, we
get \bee
A\| C^{\frac{1}{2}}f\|^2 &\leq  A\lambda^{'} \|K^*C^{\frac{1}{2}}f\|^2\\
& \leq \lambda^{'}  \|(CS)^{\frac{1}{2}}f\|^2, \eee which implies
\bea\label{e2.2} &&\frac{A}{\lambda^{'}}\| C^{\frac{1}{2}}f\|^2
\leq
\|S^{\frac{1}{2}}C^{\frac{1}{2}}f\|^2 \nn\\
&\Rightarrow& \sqrt{\frac{A}{\lambda^{'}}}\| C^{\frac{1}{2}}f\|
\leq \|S^{\frac{1}{2}}C^{\frac{1}{2}}f\|. \eea Now by using Lemma
\ref{l3.1}, we have \bee && \langle
S^{\frac{1}{2}}C^{\frac{1}{2}}f, S^{\frac{1}{2}}C^{\frac{1}{2}}f
\rangle \geq \sqrt{\frac{A}{\lambda^{'}}} \langle
C^{\frac{1}{2}}f, C^{\frac{1}{2}}f  \rangle \\&\Rightarrow&
\langle C^{\frac{1}{2}}f, C^{\frac{1}{2}}f  \rangle \leq
\sqrt{\frac{\lambda^{'}}{A}}\langle S_Cf,f\rangle. \eee Also \bee
\langle C^\frac{1}{2}K^*f, C^\frac{1}{2}K^*f \rangle & \leq \|K^*\|^2 \langle C^\frac{1}{2}f,C^\frac{1}{2}f \rangle ~~~~~~~~\\
&\leq \|K^*\|^2  \sqrt{\frac{\lambda^{'}}{A}} \langle S_Cf,f
\rangle.~~~~~~ \eee This implies that \bea \label{e2.3}
\frac{1}{\|K^*\|^2}\sqrt{\frac{A}{\lambda^{'}}} \langle
C^\frac{1}{2}K^*f,C^\frac{1}{2}K^*f \rangle \leq \langle S_Cf,f
\rangle .\eea Since $S_C$ is positive, self adjoint and bounded
$\mathcal{A}$-linear map, we can write \bee \langle
S_C^{\frac{1}{2}}f,  S_C^{\frac{1}{2}}f\rangle = \langle S_Cf,f
\rangle = \sum_{j \in \mathbb{J}} \langle f,\psi_j \rangle \langle
C\psi_j,f \rangle, \eee and hence by using Lemma \ref{l3.2}, there
exists some $B'>0$ such that \bee \langle
S_C^{\frac{1}{2}}f,S_C^{\frac{1}{2}}f \rangle \leq B^{'} \langle
f,f \rangle \eee \bea \label{e2.4} \implies \langle S_Cf,f \rangle
\leq B^{'} \langle f,f \rangle, \forall f \in \mathcal{H}. \eea
Therefore from (\ref{e2.3}) and (\ref{e2.4}), we conclude that
$\{\psi_j\}_{j \in \mathbb{J}}$ is a $C$-controlled $K$-frame in
Hilbert $C^*$-module $\mathcal{H}$ with frame bounds
$\frac{1}{\|K^*\|^2}\sqrt{\frac{A}{\lambda^{'}}}$ and $B^{'}$.
\end{proof}

\begin{l1} \label{l3.3} \rm Let $C \in GL^{+}(\mathcal{H})$, $CS_{C} = S_{C}C$ and $R(S_{C}^{\frac{1}{2}}) \subseteq R((CS_C)^{\frac{1}{2}})$ with
 $\overline{R((S_{C}^{\frac{1}{2}})^*) }$ is orthogonally complemented. Then $\| S_C^{\frac{1}{2}}f\|^2 \leq  \lambda \|(CS_C)^{\frac{1}{2}}f\|^2$ for some $\lambda > 0$.
\end{l1}
\begin{proof} By the assumption that $R(S_C^{\frac{1}{2}}) \subseteq R((CS_C)^{\frac{1}{2}})$ with $\overline{R((S_C^{\frac{1}{2}})^*) }$ orthogonally complemented.
Then by using Theorem \ref{t3.1}, there exists some $\lambda > 0$
such that \bee (S_C^{\frac{1}{2}}){(S_C^{\frac{1}{2}})^*} \leq
\lambda ((CS_C)^{\frac{1}{2}}) ((CS_C)^{\frac{1}{2}})^*. \eee This
implies that \bee &&\big\langle
(S_C^{\frac{1}{2}}){(S_C^{\frac{1}{2}})^*}f,f \big\rangle \leq
\lambda \big\langle ((CS_C)^{\frac{1}{2}})
((CS_C)^{\frac{1}{2}})^*f,f\big\rangle\\
&\Rightarrow& \| S_C^{\frac{1}{2}}f\|^2 \leq  \lambda
\|(CS_C)^{\frac{1}{2}}f\|^2,~~\forall f\in \mathcal{H}.
 \eee
\end{proof}
In the following theorem, we prove a characterization of
$C$-controlled Bessel sequence.
\begin{t1}\rm Let $\{\psi_j\}_{j \in \mathbb{J}}$ be a sequence of a finitely or countably generated Hilbert $\mathcal{A}$-module $\mathcal{H}$ over a unital $C^*$-algebra
$\mathcal{A}$. Suppose that $C$ commutes with the controlled frame
operator $S_{C}$ and $R(S_C^{\frac{1}{2}}) \subseteq
R((CS_C)^{\frac{1}{2}})$ with $\overline{R((S_C^{\frac{1}{2}})^*)
}$ is orthogonally complemented. Then $\{\psi_j\}_{j \in
\mathbb{J}}$ is a $C$-controlled Bessel sequence with  bound $B$
if and only if the operator $U \colon l^2 \mathcal{(A)} \to
\mathcal{H} $ defined by
    \bee
    U\{a_j\}_{j \in \mathbb{J}} = \sum_{j \in \mathbb{J}} a_j C \psi_j
    \eee
    is a well defined bounded operator from $l^2 \mathcal{(A)}$ into $\mathcal{H} $ with $\|U\| \leq \sqrt{B} \|C^{\frac{1}{2}}\|$.
\end{t1}
\begin{proof} Suppose that $\{\psi_j\}_{j \in \mathbb{J}}$ is a $C$-controlled Bessel sequence with bound $B$. Therefore we have
\bee
    \|\sum_{j \in \mathbb{J}} \langle f,\psi_j \rangle  \langle C\psi_j,f \rangle \| = \|\langle S_Cf,f \rangle \| \leq B \| f \|^2, \forall f \in
    \mathcal{H}.
\eee
We first show that $U$ is a well-defined operator. For arbitrary $n>m$, we have
\bee
\|\sum_{j = 1}^{n}a_j C \psi_j - \sum_{j = 1}^{m}a_j C \psi_j\|^2 &&= \|\sum_{j = m+1}^{n}a_j C \psi_j \|^2 ~~~~~~~~~~\nonumber \\
&& = \sup_{\|f\| = 1} \big\|\big\langle \sum_{j = m+1}^{n}a_j C \psi_j, f \big\rangle\big\|^2 \nonumber \\
&& = \sup_{\|f\| = 1} \big\|\sum_{j = m+1}^{n}a_j \langle C \psi_j, f \rangle\big\|^2 \nonumber \\
&&\leq \sup_{\|f\| = 1} \big\|\sum_{j = m+1}^{n} \langle f,C\psi_j
\rangle  \langle C\psi_j,f \rangle\big\| \big\|\sum_{j =
m+1}^{n}a_ja_{j}^{*}\big\| \nonumber \eee \bee
&&~~~~~~~~~~~~~~~~~~~~~~~~~~~~~= \sup_{\|f\| = 1} \big\|\big\langle \sum_{j = m+1}^{n} \langle f,C\psi_j \rangle C\psi_j,f \big\rangle\big\| \big\|\sum_{j = m+1}^{n}a_ja_{j}^{*}\big\|  \nonumber \\
&&~~~~~~~~~~~~~~~~~~~~~~~~~~~~~\leq \sup_{\|f\| = 1} \big\|\langle CS_Cf,f \rangle \big\|\big\|\sum_{j = m+1}^{n}a_ja_{j}^{*}\big\| \nonumber \\
&&~~~~~~~~~~~~~~~~~~~~~~~~~~~~~~= \sup_{\|f\| = 1} \big\|\langle (CS_C)^{\frac{1}{2}}f, (CS_C)^{\frac{1}{2}}f \rangle \big\|\big\|\sum_{j = m+1}^{n}a_ja_{j}^{*}\big\| \nonumber \\
&&~~~~~~~~~~~~~~~~~~~~~~~~~~~~~~\leq \sup_{\|f\| = 1} \|(CS_C)^{\frac{1}{2}}f\|^2\|a_j\|^2 \nonumber \\
&&~~~~~~~~~~~~~~~~~~~~~~~~~~~~~\leq \sup_{\|f\| = 1} \|C^{\frac{1}{2}}\|^2\|S_C^{\frac{1}{2}}f\|^2\|a_j\|^2 \nonumber \\
&& ~~~~~~~~~~~~~~~~~~~~~~~~~~~~~\leq \sup_{\|f\| =
1}B\|f\|^{2}\|C^{\frac{1}{2}}\|^2 \|a_j\|^2
=B\|C^{\frac{1}{2}}\|^2 \|a_j\|^2\nonumber \eee
This shows that $\displaystyle\sum_{j \in J}a_j C \psi_j$ is a Cauchy sequence which is convergent in $\mathcal{H}$. Thus $U(\{a_j\}_{j \in \mathbb{J}})$ is
a well defined operator from $l^2 \mathcal{(A)}$ into $\mathcal{H}$.\\
For boundedness of $U$, we consider \bee
\|U\{a_j\}_{j \in \mathbb{J}}\|^2 && = \sup_{\|f\| = 1} \|\langle U\{a_j\}, f \rangle\|^2 \nonumber \\
&& = \sup_{\|f\| = 1} \big\|\sum_{j \in \mathbb{J}}a_j \langle  C \psi_j, f \rangle\big\|^2 \nonumber \\
&&\leq \sup_{\|f\| = 1} \big\|\sum_{j \in \mathbb{J}} \langle f,C\psi_j \rangle  \langle C\psi_j,f \rangle\big\| \big\|\sum_{j \in \mathbb{J}}a_ja_{j}^{*}\big\| \nonumber \\
&&= \sup_{\|f\| = 1} \big\|\big\langle \sum_{j \in \mathbb{J}} \langle f,C\psi_j \rangle C\psi_j,f \big\rangle\big\| \big\|\sum_{j \in \mathbb{J}}a_ja_{j}^{*}\big\| \nonumber \\
&&= \sup_{\|f\| = 1} \big\|\langle CS_Cf,f \rangle \big\|\big\|\sum_{j \in \mathbb{J}}a_ja_{j}^{*}\big\| \nonumber \\
&&= \sup_{\|f\| = 1} \big\|\big\langle (CS_C)^{\frac{1}{2}}f, (CS_C)^{\frac{1}{2}}f \big\rangle \big\|\big\|\sum_{j \in \mathbb{J}}a_ja_{j}^{*}\big\| \nonumber \\
&&= \sup_{\|f\| = 1} \big\|(CS_C)^{\frac{1}{2}}f\big\|^2\|a_j\|^2 \nonumber \\
&& \leq \sup_{\|f\| = 1} \|C^{\frac{1}{2}}\|^2\|S_C^{\frac{1}{2}}f\|^2\|a_j\|^2 \nonumber \\
&& \leq B \|C^{\frac{1}{2}}\|^2 \|a_j\|^2.\eee This implies that $
    \|U\| \leq \sqrt{B} \|C^{\frac{1}{2}}\|.$\\

\noindent Now assume that $U$ is well defined operator from $l^2
\mathcal{(A)}$ into $\mathcal{H} $  and $\|U\| \leq \sqrt{B}
\|C^{\frac{1}{2}}\|$. We now prove that
$\{\psi_j\}_{j \in \mathbb{J}}$ is a $C$-controlled Bessel sequence with Bessel bound $B$.\\
For arbitrary $f \in \mathcal{H} $ and $\{a_j\} \in l^2 \mathcal{(A)}$, we have
\bee
\big\langle f, U\{a_j\}\big\rangle &=  \big\langle f, \sum_{j \in \mathbb{J}} a_j C \psi_j \big\rangle\\
& = \big\langle \sum_{j \in \mathbb{J}} a_j^{*} Cf, \psi_j \big\rangle\\
& =  \sum_{j \in \mathbb{J}} \langle  Cf, \psi_j \rangle a_j^{*}.
\eee Therefore we get \bee \big\langle f, U\{a_j\}\big\rangle & =
\big\langle \{ \langle Cf, \psi_j \rangle\}, \{a_j\}\big\rangle.
\eee
This implies that $U$ is has an adjoint, and  $U^*f = \{ \langle Cf, \psi_j \rangle\}$. Also, $\|U\| = \|U^*\|$.\\
So we have \bea\label{eqn6}\|U^*f\|^2 = \|\langle U^*f,
U^*f\rangle\|=\|\langle UU^*f, f\rangle\| = \|\langle
CS_Cf,f\rangle\|& =\| (CS_C)^{\frac{1}{2}}f\|^2\nn\\  & \leq B
\|C^{\frac{1}{2}}\|^2\|f\|^2. \eea By using Lemma \ref{l3.3}, we
have $\| S_C^{\frac{1}{2}}f\|^2 \leq  \lambda
\|(CS_C)^{\frac{1}{2}}f\|^2$ for some $\lambda > 0$. Using
(\ref{eqn6}) we get \bee \| S_C^{\frac{1}{2}}f\|^2 & \leq  \lambda
\|(CS_C)^{\frac{1}{2}}f\|^2  \leq \lambda B
\|C^{\frac{1}{2}}\|^2\|f\|^2. \eee Therefore $\{\psi_j\}_{j \in
\mathbb{J}}$ is a $C$-controlled Bessel sequence with Bessel bound
$\lambda B \|C^{\frac{1}{2}}\|^2$.
\end{proof}

\begin{p1}\rm Let $\{\psi_j\}_{j \in \mathbb{J}}$ be a $C$-controlled $K$-frame in  $\mathcal{H}$. Then $ ACKK^{*}I  \leq  S_{c} \leq BI$.
\end{p1}
\begin{proof} Suppose $\{\psi_j\}_{j \in \mathbb{J}}$  is a $C$-controlled $K$-frame with bounds $A$ and $B$. Then
\bee &&A\langle C^{\frac{1}{2}}K^{*}f, C^{\frac{1}{2}}K^{*}f
\rangle \leq \sum_{j \in \mathbb{J}} \langle f,\psi_j \rangle
\langle C\psi_j,f \rangle \leq B \langle f,f \rangle , \forall f
\in
\mathcal{H}.\\
&\Rightarrow&  A\langle CKK^{*}f, f \rangle  \leq \langle S_{C}f,f
\rangle \leq B \langle f,f \rangle.\\
&\Rightarrow& ACKK^{*}I  \leq S_{C} \leq BI. \eee
\end{proof}

\begin{p1}\rm Let $\{\psi_j\}_{j \in \mathbb{J}}$ be a $C$-controlled Bessel sequence in $\mathcal{H}$ and $C \in GL^{+}(\mathcal{H})$.
Then $\{\psi_j\}_{j \in \mathbb{J}}$ is a $C$-controlled $K$-frame for $\mathcal{H}$, if and only if there exists $A > 0$ such
that $CS \geq ACKK^*$.
\end{p1}
\begin{proof} The sequence $\{\psi_j\}_{j \in \mathbb{J}}$ is a controlled $K$-frame for $\mathcal{H}$ with frame bounds $A$, $B$ and frame operator $S_C$, if and only if
\bee
&& A\langle C^{\frac{1}{2}}K^{*}f, C^{\frac{1}{2}}K^{*}f \rangle  \leq \sum_{j \in \mathbb{J}} \langle f,\psi_j \rangle  \langle C\psi_j,f \rangle \leq B \langle f,f \rangle ,
\forall f \in \mathcal{H}.\\
&& \Leftrightarrow  A\langle CKK^{*}f, f \rangle  \leq \langle S_Cf,f \rangle \leq B \langle f,f \rangle. \\
&& \Leftrightarrow  A\langle CKK^{*}f, f \rangle  \leq \langle CSf,f \rangle \leq B \langle f,f \rangle. \\
&&\Leftrightarrow  ACKK^{*}I \leq CS. \eee
\end{proof}

In the following two propositions we establish the
inter-relationship between $K$-frame and $C$-controlled $K$-frame.

\begin{p1}\rm Let $C \in GL^{+}(\mathcal{H})$, $K\in L(\mathcal{H})$, $KC=CK$, $R(C^{\frac{1}{2}}) \subseteq R(K^*C^{\frac{1}{2}})$ with $\overline{R((C^{\frac{1}{2}})^*) }$
is orthogonally complemented, and $\{\psi_j\}_{j \in \mathbb{J}}$ be a $C$-controlled $K$-frame for
$\mathcal{H}$ with lower and upper frame bounds $A$ and $B$,
respectively. Then $\{\psi_j\}_{j \in \mathbb{J}}$ is a $K$-frame
for $\mathcal{H}$ with lower and upper frame bounds
$A\|C^{\frac{1}{2}}\|^{-2}$ and $B\|C^{\frac{-1}{2}}\|^2$,
respectively.
\end{p1}
\begin{proof} Suppose $\{\psi_j\}_{j \in \mathbb{J}}$ is a $C$-controlled $K$-frame for $\mathcal{H}$ with bound $A$ and $B$. Then by Theorem \ref{t3.2}, we have
\bee
    A\| C^{\frac{1}{2}}K^{*}f\|^2  \leq \|\sum_{j \in \mathbb{J}} \langle f,\psi_j \rangle  \langle C\psi_j,f \rangle \| \leq B \| f \|^2,  \forall f \in \mathcal{H}.
\eee
Now,
\begin{equation}
\begin{split}
A\|K^*f\|^2 & =  A\|C^{\frac{-1}{2}}C^{\frac{1}{2}}K^*f\|^2 \nonumber \\
& \leq A\|C^{\frac{1}{2}}\|^2\|C^{\frac{-1}{2}}K^*f\|^2 \nonumber \\
& \leq \|C^{\frac{1}{2}}\|^2\|\sum_{j \in \mathbb{J}} \langle
f,\psi_j\rangle \langle \psi_j,f\rangle\|.
\end{split}
\end{equation}
This implies that
\bee
A\|C^{\frac{1}{2}}\|^{-2}\|K^*f\|^2 \leq
\|\sum_{j \in \mathbb{J}} \langle f,\psi_j\rangle \langle \psi_j,f
\rangle\|
\eee
On the other hand for every $f \in
\mathcal{H}$, \bee
\|\sum_{j \in \mathbb{J}} \langle f,\psi_j\rangle \langle \psi_j,f \rangle\|  &&= \|\langle Sf,f\rangle\| \nonumber \\
 &&=\|\langle C^{-1}CSf,f\rangle\| \nonumber \\
&& =\|\langle(C^{-1}CS)^{\frac{1}{2}}f,(C^{-1}CS)^{\frac{1}{2}}f\rangle\| \nonumber \\
&& =\|(C^{-1}CS)^{\frac{1}{2}}f\|^2 \nonumber \\
&& \leq\|C^{\frac{-1}{2}}\|^2\|(CS)^{\frac{1}{2}}f\|^2 \nonumber \\
&& = \|C^{\frac{-1}{2}}\|^2\|\langle(CS)^{\frac{1}{2}}f,(CS)^{\frac{1}{2}}f \rangle\| \nonumber \\
&& = \|C^{\frac{-1}{2}}\|^2\|\langle CSf,f \rangle\| \nonumber \\
&& \leq \|C^{\frac{-1}{2}}\|^2B\|f\|^2. \nonumber \eee Therefore,
$\{\psi_j\}_{j \in \mathbb{J}}$ is a $K$-frame with lower and
upper frame bounds $A\|C^{\frac{1}{2}}\|^{-2}$ and
$B\|C^{\frac{-1}{2}}\|^2$, respectively.
\end{proof}

\begin{p1}\rm Let  $C \in GL^{+}(\mathcal{H})$, $K \in L(\mathcal{H})$, $KC=CK$, $R(C^{\frac{1}{2}}) \subseteq R(K^*C^{\frac{1}{2}})$ with $\overline{R((C^{\frac{1}{2}})^*) }$
is orthogonally complemented. Let  $\{\psi_j\}_{j \in \mathbb{J}}$ be a $K$-frame
    for $\mathcal{H}$ with lower and upper frame bounds $A$ and $B$,
    respectively. Then $\{\psi_j\}_{j \in \mathbb{J}}$ is a
    $C$-controlled $K$-frame for $\mathcal{H}$ with lower and upper
    frame bounds $A$ and $\|C\|\|S\|$, respectively.
\end{p1}
\begin{proof} Suppose $\{\psi_j\}_{j \in \mathbb{J}}$ is a $K$-frame with frame bounds $A$ and $B$. Then by equivalence condition \cite{Gvruta} of K-frame, we have
    \bee
    A\|K^{*}f\|^2  \leq \|\sum_{j \in \mathbb{J}} \langle f,\psi_j \rangle  \langle \psi_j,f \rangle \| \leq B \| f \|^2, \forall f \in \mathcal{H}.
    \eee For any $f\in \mathcal{H}$,
    \bea\label{e2.5}
    A\|C^{\frac{1}{2}}K^*f\|^2 &=& A\|K^*C^{\frac{1}{2}}f\|^2\nn  \\
    &\leq&\|\sum_{j \in \mathbb{J}} \langle
    C^{\frac{1}{2}}f,\psi_j\rangle \langle \psi_j,C^{\frac{1}{2}}f
    \rangle\| \nn\\
    &=& \|\sum_{j \in \mathbb{J}} \langle
    C^{\frac{1}{2}}f,\psi_j\rangle \psi_j,C^{\frac{1}{2}}f
    \rangle\|\nn
    \\
    &=& \|\langle C^{\frac{1}{2}}Sf, C^{\frac{1}{2}}f \rangle\|\nn \\
    &=& \|\langle CSf,f \rangle\|. \eea
    On the other hand for every $f \in \mathcal{H}$,
    \bea\label{e2.6}
    \|\langle CSf,f \rangle\| &=& \|\langle Sf,C^*f \rangle\| \nn\\
    &=& \|\langle Sf,Cf \rangle\|\nn \\
    &\leq&  \| Sf\|\|Cf\|\nn\\
    &\leq&  \|C\| \|S\|\|f\|^2. \eea
    Therefore from (\ref{e2.5}),(\ref{e2.6}) and Theorem \ref{t3.2}, we conclude that
    $\{\psi_j\}_{j \in \mathbb{J}}$ is a $C$-controlled $K$-frame with
    bounds $A$ and $\|C\|\|S\|$.
\end{proof}

\begin{t1}\rm Let $C \in GL^{+}(\mathcal{H})$, $\{\psi_j\}_{j \in \mathbb{J}}$ be a $C$-controlled $K$-frame for $\mathcal{H}$ with bounds $A$ and $B$.
Let $M, K \in L(\mathcal{H})$ with $R(M) \subset R(K)$,
$\overline{R(K^*)}$ orthogonally complemented, and $C$ commutes
with $M$ and $K$ both. Then $\{\psi_j\}_{j \in \mathbb{J}}$ is a
$C$-controlled $M$-frame for $\mathcal{H}$.
\end{t1}
\begin{proof} Suppose $\{\psi_j\}_{j \in \mathbb{J}}$ is a $C$-controlled $K$-frame for $\mathcal{H}$ with bounds $A$ and $B$. Then
\bea\label{eqn7}
    A\langle C^{\frac{1}{2}}K^{*}f, C^{\frac{1}{2}}K^{*}f \rangle  \leq \sum_{j \in \mathbb{J}} \langle f,\psi_j \rangle  \langle C\psi_j,f \rangle \leq B \langle f,f \rangle,
     ~~\forall f \in \mathcal{H}.
\eea Since $R(M) \subset R(K)$, from Theorem \ref{t3.1}, there
exists some $\lambda^{'} > 0$ such that $MM^* \leq
\lambda^{'}KK^*$. So we have
\begin{equation}
\langle MM^*C^{\frac{1}{2}}f,C^{\frac{1}{2}}f \rangle \leq \lambda^{'} \langle KK^*C^{\frac{1}{2}}f,C^{\frac{1}{2}}f \rangle . \nonumber \\
\end{equation}
Multiplying the above inequality by $A$, we get
\begin{center}
    $\frac{A}{\lambda^{'}} \langle MM^*C^{\frac{1}{2}}f,C^{\frac{1}{2}}f \rangle \leq A \langle KK^*C^{\frac{1}{2}}f,C^{\frac{1}{2}}f \rangle. $
\end{center}
From (\ref{eqn7}), we have
\begin{center}
    $\frac{A}{\lambda^{'}} \langle MM^*C^{\frac{1}{2}}f,C^{\frac{1}{2}}f \rangle \leq \sum_{j \in \mathbb{J}} \langle f,\psi_j \rangle
    \langle C\psi_j,f \rangle \leq B \langle f,f \rangle$, for all $f \in
    \mathcal{H}$.
\end{center}
Therefore, $\{\psi_j\}_{j \in \mathbb{J}}$ is a $C$-controlled
$M$-frame  with lower and upper frame bounds
$\frac{A}{\lambda^{'}} $ and $B$, respectively.
\end{proof}

In the following result, we investigate the invariance of a
$C$-controlled Bessel sequence under a adjointable operator.

\begin{p1}\label{p1}\rm Let $\{\psi_j\}_{j \in \mathbb{J}}$ be a $C$-controlled Bessel sequence with bound $D$. Let $T \in L(\mathcal{H})$ and $CT = TC$.
Then $\{T\psi_j\}_{j \in \mathbb{J}}$ is also $C$-controlled
Bessel sequence with bound  $D {\|T^*\|}^{2}$.
\end{p1}
\begin{proof} Suppose $\{\psi_j\}_{j \in \mathbb{J}}$ is a $C$-controlled Bessel sequence with bound $D$. Then we have
\bee
    \sum_{j \in \mathbb{J}} \langle f,\psi_j \rangle  \langle C\psi_j,f \rangle \leq D \langle f,f \rangle, \forall f \in \mathcal{H}.
\eee For every $f\in \mathcal{H}$,
\begin{equation}
\begin{split}
\sum_{j \in \mathbb{J}} \langle f, T\psi_j \rangle \langle CT\psi_j, f \rangle & = \sum_{j \in \mathbb{J}} \langle T^*f, \psi_j \rangle \langle TC\psi_j, f \rangle \\
& = \sum_{j \in \mathbb{J}} \langle T^*f, \psi_j \rangle \langle C\psi_j, T^*f \rangle  \nonumber \\
& \leq D \langle T^*f, T^*f \rangle  \nonumber \\
& \leq D {\|T^*\|}^{2} \langle f,f \rangle.  \nonumber \\
\end{split}
\end{equation}
Thus $\{T\psi_j\}_{j \in \mathbb{J}}$ is also $C$-controlled
Bessel sequence with bound  $D {\|T^*\|}^{2}$.
\end{proof}

Now, we investigate the invariance of a $C$-controlled $K$-frame
under a adjointable operator.

\begin{t1}\rm Let  $C \in GL^{+}(\mathcal{H})$, $K \in L(\mathcal{H})$ and $\{\psi_j\}_{j \in \mathbb{J}}$ be a $C$-controlled $K$-frame for $\mathcal{H}$
with lower and upper bounds $A$ and $B$, respectively. If $T \in
L(\mathcal{H})$ with closed range such that $\overline{R(TK)}$ is
orthogonally complemented and $C, K, T$ commute with each other.
Then $\{T\psi_j\}_{j \in \mathbb{J}}$ is a $C$-controlled
$K$-frame for $R(T)$.
\end{t1}
\begin{proof} Suppose $\{\psi_j\}_{j \in \mathbb{J}}$ is a $C$-controlled $K$-frame for $\mathcal{H}$ with bound $A$ and $B$. Then
\bee
    A\langle C^{\frac{1}{2}}K^{*}f, C^{\frac{1}{2}}K^{*}f \rangle  \leq \sum_{j \in \mathbb{J}} \langle f,\psi_j \rangle  \langle C\psi_j,f \rangle \leq B \langle f,f \rangle,
    \forall f \in \mathcal{H}.
\eee We know that if $T$ has closed range then $T$ has
Moore-Penrose inverse $T^{\dagger}$ such that $TT^{\dagger}T=T$
and $T^{\dagger}TT^{\dagger}=T^{\dagger}$. So
$TT^{\dagger}|_{R(T)}=I_{R(T)}$ and
$(TT^{\dagger})^*=I^*=I=TT^{\dagger}.$\\
We have
\begin{equation}
\begin{split}
\langle K^*C^{\frac{1}{2}}f, K^*C^{\frac{1}{2}}f \rangle &= \big\langle (TT^{\dagger})^{*} K^*C^{\frac{1}{2}}f, (TT^{{\dagger}})^{*} K^*C^{\frac{1}{2}}f \big\rangle \nonumber \\
&= \big\langle T^{{\dagger}*}T^*K^*C^{\frac{1}{2}}f, T^{{\dagger}*}T^* K^*C^{\frac{1}{2}}f \big\rangle \nonumber \\
&\leq \|(T^{\dagger})^*\|^{2} \big\langle T^*K^*C^{\frac{1}{2}}f, T^* K^*C^{\frac{1}{2}}f \big\rangle. \nonumber \\
\end{split}
\end{equation}
This  implies that \bea\label{eqn8}
    \|(T^{\dagger})^*\|^{-2} \big\langle K^*C^{\frac{1}{2}}f, K^*C^{\frac{1}{2}}f \big\rangle \leq \big\langle T^*K^*C^{\frac{1}{2}}f, T^* K^*C^{\frac{1}{2}}f
    \big\rangle.
\eea Since $R(T^*K^*) \subset R(K^*T^*)$, by using Theorem
\ref{t3.1}, there exists some $\lambda^{'} > 0$ such that
\bea\label{eqn9}\big\langle
T^*K^*C^{\frac{1}{2}}f,T^*K^*C^{\frac{1}{2}}f \big\rangle \leq
\lambda^{'} \big\langle
K^*T^*C^{\frac{1}{2}}f,K^*T^*C^{\frac{1}{2}}f \big\rangle. \eea
Therefore, using (\ref{eqn8}) and (\ref{eqn9}) we get
\begin{equation}
\begin{split}
\sum_{j \in \mathbb{J}} \langle f,T\psi_j \rangle  \langle
CT\psi_j,f \rangle & = \sum_{j \in \mathbb{J}} \langle T^*f,\psi_j
\rangle  \langle TC\psi_j,f \rangle \\ \nonumber & =  \sum_{j \in
\mathbb{J}} \langle T^*f,\psi_j \rangle  \langle C\psi_j,T^*f
\rangle \\ \nonumber & \geq  A \big\langle
C^{\frac{1}{2}}K^*T^*f,C^{\frac{1}{2}}K^*T^*f \big\rangle  \\
\nonumber & \geq  A \lambda^{'}\langle
T^*C^{\frac{1}{2}}K^*f,T^*C^{\frac{1}{2}}K^*f \rangle  \\
\nonumber & \geq  A \lambda^{'}\|(T^{\dagger})^*\|^{-2}\langle
C^{\frac{1}{2}}K^*f, C^{\frac{1}{2}}K^*f \rangle.  \nonumber
\end{split}
\end{equation}
This gives the lower frame inequality for $\{T\psi_j\}_{j \in
\mathbb{J}}$. On the other hand by Proposition \ref{p1},
$\{T\psi_j\}_{j \in \mathbb{J}}$ is a $C$-controlled Bessel
sequence. So $\{T\psi_j\}_{j \in \mathbb{J}}$ is a $C$-controlled
$K$-frame for $R(T)$.
\end{proof}

\begin{t1} \rm
    Let  $C \in GL^{+}(\mathcal{H})$, $K \in L(\mathcal{H})$ and $\{\psi_j\}_{j \in \mathbb{J}}$ be a $C$-controlled $K$-frame for $\mathcal{H}$ with lower and upper
     bound $A$, $B$ respectively. If $T \in L(\mathcal{H})$ is a isometry such that $R(T^{*}K^{*}) \subset R(K^{*}T^{*})$ with  $\overline{R(TK)}$ is orthogonally
     complemented and $C, K, T$ commute with each other. Then $\{T\psi_j\}_{j \in \mathbb{J}}$ is a $C$-controlled $K$-frame for $\mathcal{H}$.
\end{t1}
\begin{proof} By Theorem \ref{t3.1}, there exist some $\lambda > 0$ such that $\|T^{*}K^{*}C^{\frac{1}{2}}f\|^2 \leq \lambda \|K^{*}T^{*}C^{\frac{1}{2}}f\|^2$. Suppose
$A$ is a lower bound for the $C$-controlled $K$-frame $\{\psi_j\}_{j \in \mathbb{J}}$. Since $T$ is an isometry, then \\
\bea \label{e2.7}
\frac{A}{\lambda} \|C^{\frac{1}{2}}K^{*}f\|^2  \nn
& = \frac{A}{\lambda} \|T^{*}C^{\frac{1}{2}}K^{*}f\|^2~~~~~~~~~~~~~~\\ \nn
& \leq A\|K^{*}T^{*}C^{\frac{1}{2}}f\|^2~~~~~~~~~~~~~\\ \nn
& = A\|C^{\frac{1}{2}}K^{*}T^{*}f\|^2 ~~~~~~~~~~~~~\\ \nn
& \leq \sum_{j \in \mathbb{J}} \langle T^*f,\psi_j
\rangle  \langle C\psi_j,T^{*}f \rangle \\ \nn
& = \sum_{j \in \mathbb{J}} \langle f,T\psi_j
\rangle  \langle TC\psi_j,f \rangle ~~~~\\
& = \sum_{j \in \mathbb{J}} \langle f,T\psi_j \rangle  \langle
CT\psi_j,f \rangle~~~~ \eea Therefore from Proposition \ref{p1}
and inequality (\ref{e2.7}), we conclude that $\{T\psi_j\}_{j \in
\mathbb{J}}$ is a $C$-controlled $K$-frame for $\mathcal{H}$ with
bounds $\frac{A}{\lambda}$ and $B {\|T^*\|}^{2}$.
    \end{proof}

Now we prove a perturbation result for $C$-controlled $K$-frame.

\begin{t1}\rm Let $F = \{f_j\}_{j \in \mathbb{J}}$ be a $C$-controlled $K$-frame for $\mathcal{H}$ , with controlled frame operator $S_C$. Suppose $K \in L(\mathcal{H})$, $KC=CK$, $R(C^{\frac{1}{2}}) \subseteq R(K^*C^{\frac{1}{2}})$ with $\overline{R((C^{\frac{1}{2}})^*) }$ is orthogonally complemented .
If $G = \{g_j\}_{j \in \mathbb{J}}$ is a non zero sequence in
$\mathcal{H}$, and $E = T_F - T_G$ be a compact operator, where
$T_G(\{c_j\}_{j \in \mathbb{J}}) = \sum_{j\in \mathbb{J}}c_{j}
g_{j}$ for $\{c_j\}_{j \in \mathbb{J}} \in l^2 \mathcal{(A)}$,
then $G = \{g_j\}_{j \in \mathbb{J}}$ is a $C$-controlled
$K$-frame for $\mathcal{H}$.
\end{t1}
\begin{proof} Let $\{f_j\}_{j \in \mathbb{J}}$ be a $C$-controlled $K$-frame with bounds  $A$ and $B$, then because of \\Theorem \ref{t3.2}, we have
\bee
    A\| C^{\frac{1}{2}}K^{*}f\|^2  \leq \|\sum_{j \in \mathbb{J}} \langle f,f_j \rangle  \langle C f_j,f \rangle \| \leq B \| f \|^2, \forall f \in
    \mathcal{H}.
\eee
This implies $\|T_F\|^2 \leq B\|C^{\frac{-1}{2}}\|^2$.\\
Let $V = T_F - E$ be an operator from $l^2 \mathcal{(A)}$ into $\mathcal{H}$. Since $T_F$ and $E$ are bounded, then the operator $V$ is bounded. Therefore $\|V\|=\|V^*\|$.\\
For any $f \in \mathcal{H}$,
\begin{equation}
\begin{split}
V^*f & = T_F^*f - E^*f \\
& = \{\langle f, f_j \rangle\}_{j \in \mathbb{J}} - \{\langle f, f_j - g_j \rangle\}_{j \in \mathbb{J}}\\
& = \{\langle f, f_j \rangle\}_{j \in \mathbb{J}} - \{\langle f_j - g_j, f \rangle^{*}\}_{j \in \mathbb{J}}\\
& = \{\langle f, f_j \rangle\}_{j \in \mathbb{J}} - \{\langle f_j, f \rangle ^{*} - \langle g_j, f \rangle^{*} \}_{j \in \mathbb{J}}\\
& = \{\langle f, f_j \rangle\}_{j \in \mathbb{J}} - \{\langle f, f_j \rangle - \langle f, g_j \rangle \}_{j \in \mathbb{J}}\\
& = \{\langle f, g_j \rangle\}_{j \in \mathbb{J}}. \nonumber
\end{split}
\end{equation}
We have \bea\label{eqn10} V(\{c_j\}_{j \in \mathbb{J}}) = \sum_{j
\in \mathbb{J}} c_{j} g_{j}, {~\rm ~and~}S_G = VV^* .\eea

Now using (\ref{eqn10}), we have \bea\label{eqn11} \|\langle
f,CS_Gf\rangle\| = \|\langle f,CVV^*f\rangle\|
&=& \|\langle C^{\frac{1}{2}}V f,C^{\frac{1}{2}}Vf\rangle\|\nn \\
&=& \| C^{\frac{1}{2}}Vf \|^2~~~~~\nn \\
&\leq& \| C^{\frac{1}{2}}\|^2 \| Vf \|^2\nn\\
&=& \| C^{\frac{1}{2}}\|^2 \| (T_F - E)f \|^2\nn\\
&\leq& \| C^{\frac{1}{2}}\|^2 \| T_F - E \|^2 \|f\|^2\nn\\
&\leq& (\| T_F \|^2 + 2\| T_F \|\| E \|+\| E \|^2) \|
C^{\frac{1}{2}}\|^2  \|f\|^2\nn\\
 &\leq&
\big(B\|C^{\frac{-1}{2}}\|^2 +2 \sqrt{B}\|C^{\frac{-1}{2}}\|\| E
\|+\| E \|^2\big) \| C^{\frac{1}{2}}\|^2  \|f\|^2\nn\\ &=& B
\Big(\| C^{\frac{-1}{2}}\| + \frac{\|E\|}{\sqrt{B}}\Big)^2 \|
C^{\frac{1}{2}}\|^2  \|f\|^2.  \eea This inequality shows that
$\{g_j\}_{j \in \mathbb{J}}$ is a controlled Bessel sequence with
bound\\ $B \Big(\| C^{\frac{-1}{2}}\| +
\frac{\|E\|}{\sqrt{B}}\Big)^2 \|
C^{\frac{1}{2}}\|^2$.\\
Again we have \bee
VV^* &= (T_{F} - E)(T_{F} - E)^* ~~~~~~~~~~~~~~~~~~~~~~~~~~~\\
&= (T_{F} - E)(T_{F}^* - E^*) ~~~~~~~~~~~~~~~~~~~~~~~~~~~\\
&= T_{F}T_{F}^* - T_{F}E^* - ET_{F}^* + EE^*~~~~~~~~~~~~~~~\\
&= S_{F} - T_{F}E^* - ET_{F}^* + EE^*~~~~~~~~~~~~~~~~~~~
\eee

\noindent Since $E, T_{F}$ and $S_F$ are compact operators, then
$S_{F} - T_{F}E^* - ET_{F}^* + EE^*$ is a compact operator.
Therefore $S_{F} - T_{F}E^* - ET_{F}^* + EE^* + I$ is a bounded
operator with closed range. Thus, $ VV^* = S_{F} - T_{F}E^* -
ET_{F}^* + EE^*$ is a bounded operator with closed range. Also
$VV^{*}$ is injective as $V$ is injective. Hence $VV^{*}(=S_{G})$
is bounded below. So there exists some constant $A > 0$ such that
\bea \label{eqn13} A\|C^{\frac{1}{2}}f\| \leq
\|S_{G}C^{\frac{1}{2}}f\|. \eea Now \bee\label{eqn14}
\|C^{\frac{1}{2}}K^*f\|^2 &=& \|K^*C^{\frac{1}{2}}f\|^2\nn\\
&\leq& \|K^*\|^2 \|C^{\frac{1}{2}}f\|^2\nn\\
&\leq& \frac{1}{A^2}\|K^*\|^2\|S_{G}C^{\frac{1}{2}}f\|^2. \eee
This implies that \bea\label{eqn15}
\frac{A^2}{\|K^*\|^2}\|C^{\frac{1}{2}}K^*f\|^2 \leq
\|S_{G}C^{\frac{1}{2}}f\|^2. \eea Therefore from (\ref{eqn11}) and
(\ref{eqn15}), we conclude that $G = \{g_j\}_{j \in \mathbb{J}}$
is a $C$-controlled $K$-frame for $\mathcal{H}$ with frame bounds
$\frac{A^2}{\|K^*\|^2}$ and $B \big(\| C^{\frac{-1}{2}}\| +
\frac{\|E\|}{\sqrt{B}}\big)^2 \| C^{\frac{1}{2}}\|^2$.
\end{proof}

\end{document}